\documentclass[12pt,twoside,english]{article}
\usepackage{latexsym,amsmath,amssymb}
\usepackage{babel}
\usepackage{color}
\title{Multivariate semi-discrete sampling type operators: pointwise approximation properties}
\author{Carlo Bardaro - Ilaria Mantellini}
\date{\centerline{\small  Department of Mathematics and Informatics, 
University of Perugia} \centerline{\small Via Vanvitelli 1, 06123 Perugia, 
ITALY} \centerline{\small e-mail: carlo.bardaro@unipg.it, mantell@dmi.unipg.it}}

\begin{document}
\maketitle
\vskip0,5cm

{\small{\bf Abstract:} In this paper multivariate extension of the generalized Durrmeyer sampling type series are considered. We establish a Voronovskaja type formula and a quantitative version. Finally some particular examples are discussed.}

\section{\bf Introduction}

The so-called ``generalized sampling series'' of a function $f:\mathbb{R}\rightarrow \mathbb{R}$ was introduced by P.L. Butzer and his school in Aachen during the eighties and the beginning of nineties (see e.g. \cite{RS}, \cite{BS1}, \cite{BS0}). It is defined by
\begin{eqnarray}
(S^\varphi_nf) (x):= \sum_{k=-\infty}^\infty \varphi (nx -k)f(\frac{k}{n}),\quad n \in \mathbb{N}, x \in \mathbb{R},
\end{eqnarray}
where $\varphi$ is a ``kernel'' function satisfying classical assumptions of singularity. This operator represents an extension of the classical Shannon sampling series and has a great importance in signal analysis. In particular when one considers a kernel $\varphi$ with a compact support contained in the positive real semi-axis $\mathbb{R}^+$ we can obtain efficient mathematical models for the prediction of the signal $f.$
In the space of uniformly continuous and bounded functions $C^0(\mathbb{R})$ the operator (1) behaves in a very good way: it is bounded as operator $S^\varphi_n: C^0(\mathbb{R}) \rightarrow C^0(\mathbb{R})$ and $S^\varphi_nf$ converges uniformly to $f.$

Later on generalizations to the multivariate setting were studied (see \cite{FS} and \cite{BFS}). However, when we go beyond the space $C^0(\mathbb{R}),$ by considering e.g. the Lebesgue spaces $L^p(\mathbb{R}),$ the operator $S^\varphi_n$ is no longer bounded in $L^p$ and also the convergence property (in $L^p-$sense) holds only for $f$ in suitable subspaces (\cite{BBSV1}, \cite{BMSVV}). The study of operators like (1) in $L^p-$spaces is very important in order to obtain mathematical models in signal and image processing. Therefore it appears very useful to construct operators with better approximation properties also in more general functional spaces. For this reason, taking inspiration from the theory of the Bernstein polynomials, a Kantorovich modification was introduced (see \cite{BBSV}) in which the sample values $f(k/n)$ are replaced by a mean value of $f$ in a small interval, namely we put
$$(\widetilde{S}^\varphi_nf)(x) := \sum_{k=-\infty}^\infty \varphi (nx -k)\frac{1}{n+1}\int_{k/n}^{(k+1)/n}f(u)du, \quad n \in \mathbb{N}, x \in \mathbb{R}.$$
These ``semi-discrete'' operators  are now bounded in $L^p$ (and more generally also in Orlicz spaces) and have good approximation properties. For the literature about these operators and their properties both in one and multidimensional case,  along with their applications in signal and image processing see e.g. \cite{VZ}, \cite{CV1}, \cite{CV2}, \cite{KI}. 

A more convenient version of the Bernstein polynomials was then introduced by Durrmeyer (see \cite{DU}, \cite{DE}), by replacing the integral means by a convolution integral, with the same function $\varphi$ as a kernel. It turns out that a similar modification gives operators which have better approximation properties. Therefore, in \cite{BM3}, \cite{BFM} we have introduced a similar modification for the sampling series, by using two different kernel functions, namely we defined:
$$(S^{\varphi, \psi}_nf)(x) := \sum_{k=-\infty}^\infty \varphi (nx -k)n\int_{-\infty}^\infty \psi(nu-k) f(u)du,\quad n \in  \mathbb{N}, x \in \mathbb{R},$$
where $\varphi$ and $\psi$ are two kernel functions satisfying the classical assumption of an approximate identity (see e.g. \cite{BN}) and $f$ belongs to some function spaces. In the quoted papers \cite{BM3}, \cite{BFM} we were interested in establishing certain asymptotic expressions and Voronovskaja-type formulae for the pointwise and uniform convergence within the frame of continuous functions.

Due to its importance in applications, in this paper we present a multivariate version of the operator $S^{\varphi, \psi}_n,$ acting in the space $C^0(\mathbb{R}^N)$ of all bounded and uniformly continuous functions defined on $\mathbb{R}^N$ and we study its pointwise and uniform asymptotic behaviour, along with Voronovskaja type results, in case the function $f$ satisfies some local regularity assumptions. In Section 2 we introduce the main notations and definitions and in Section 3 we give our main approximation theorems. In Section 4 we give certain quantitative estimates of the uniform convergence, employing the classical Peetre K-functional (see \cite{PE2},  \cite{DL}, \cite{AG}). Here, the main tool is an estimate of the remainder in the Taylor formula for $f$ in terms of the $K-$functional, which we extend to the multidimensional frame (Lemma 1). These results imply also direct estimates of the uniform approximation in suitable subspaces. Section 5 contains some applications of the theory developed here, to particular examples, involving Bochner-Riesz kernels (\cite{SW}) and certain multivariate box splines.

\section{Notations and definitions}

We will denote by $\mathbb{N}$~ and $\mathbb{N}_0$~the sets of positive and non negative integers respectively, by $\mathbb{Z}$~ the set of integers, by $\mathbb{R},$~ $\mathbb{R}^+,$~ $\mathbb{R}^+_0$~the sets of all real, positive real and non negative real numbers respectively.

Let $\mathbb{Z}^N,$~ $\mathbb{N}^N$~ and $\mathbb{N}^N_0$~ be the sets of all N-tuples ${\tt k}= (k_1, \ldots, k_N)$, of all integers, positive integers and non negative integers respectively  and by $\mathbb{R}^N,$ $\mathbb{R}^N_+$ the sets of all N-tuples of real numbers and positive real numbers.
Given ${\tt x}= (x_1, \ldots, x_N)$ and ${\tt y}=(y_1, \ldots, y_N)$ we will denote by 
$\|{\tt x}\|$ the euclidean norm of the vector ${\tt x},$ ${\tt x}\cdot {\tt y}$ the scalar product,
${\tt x y}=(x_1 y_1, \ldots, x_Ny_N),$ ${\tt x/y}= (x_1/y_1, \ldots, x_N/y_N),$ $1/{\tt x} = (1/x_1, \ldots, 1/x_N),$
${\tt x}^{\tt y}= \prod_{i=1}^N x_i^{y_i},$ when the power is well defined, and 
$[{\tt x}]= (|x_1|, |x_2|, \ldots, |x_N|).$

Moreover we will put by $\langle{\tt x}\rangle = \prod_{i=1}^N x_i$ and ${\tt x}\leq {\tt y}$ means that
$ x_i \leq y_i$ for every $i=1, \ldots, N.$ By ${\tt 0}$ we denote the null vector, by ${\tt 1}=(1, \ldots, 1)$ and by ${\tt e}_i,~ i=1, \ldots, N$ the vectors of the standard basis of $\mathbb{R}^N.$
Further, standard multi-index notation is used, i.e., for ${\tt k}= (k_1, \ldots, k_N) \in \mathbb{N}^N_0,$ we write 
${\tt k}! = k_1! k_2! k_N!$ and $|{\tt k}| = k_1 + \ldots + k_N.$

For a function $f: \mathbb{R}^N \rightarrow \mathbb{R},$ we denote
$$D^{{\tt k}}f:= \frac{\partial^{|{\tt k}|}}{\partial {\tt x}^{{\tt k}}}f = \frac{\partial^{|{\tt k}|}}{\partial {x_1}^{k_1}\ldots \partial {x_N}^{k_N}} f~~~(|{\tt k}| = r)$$
the r-th order derivatives of $f.$ 

Let us denote by $L^1=L^1(\mathbb{R}^N )$ the space of all Lebesgue integrable functions $f: \mathbb{R}^N \rightarrow \mathbb{R},$~ provided with the usual norm $\|f\|_1,$ by
$L^\infty = L^\infty(\mathbb{R}^N )$ the space of all the essentially bounded functions $f$~endowed with the usual supnorm $\|f\|_\infty,$ by  $C^0 = C^0(\mathbb{R}^N)$~ the subspace of all uniformly continuous and bounded functions and for $\nu\geq 1$~ by $C^\nu = C^\nu(\mathbb{R}^N)$~ the subspace of $C^0$~ whose elements $f$~ are $\nu-$times continuosly differentiable and the $\nu$-th-order derivatives are in $C^0.$ 

In what follows, we will say that a function $f$~ belongs to $C^\nu$~ locally at a point ${\tt x}\in \mathbb{R}^N$~ if there is a neighbourhood $U$~ of ${\tt x}$~ such that $f$~ is $(\nu-1)-$fold continuously differentiable in $U$~ and the $\nu$-th derivatives exist and are continuous at the point ${\tt x}.$

Let us consider two functions $\varphi \in C^0$~ and $\psi \in L^1(\mathbb{R}^N )$. 
For any $\nu\in \mathbb{N}_0$~ and ${\tt h}\in \mathbb{N}_0^N,~{\tt h}=(h_1, \ldots, h_N), $~with $|{\tt h}|= \nu$ 
let us define the algebraic moments
$$m^{\nu}_{\tt h}(\varphi, {\tt u}) := \sum_{{\tt k}\in \mathbb{Z}^N} \varphi ({\tt u}-{\tt k}) ({\tt k}-{\tt u})^{\tt h}, \quad \widetilde{m}^{\nu}_{\tt h}(\psi) := \int_{\mathbb{R}^N} {\tt t}^{\tt h} \psi({\tt t})d{\tt t}$$
and the absolute moments
$$ M^{\nu}_{\tt h}(\varphi) := \sup_{{\tt u}\in \mathbb{R}^N} \sum_{{\tt k}\in \mathbb{Z}^N} |\varphi({\tt u}-{\tt k})| [{\tt u}-{\tt k}]^{\tt h}, \quad \widetilde{M}^{\nu}_{\tt h}(\psi) := \int_{\mathbb{R}^N} [{\tt t}]^{\tt h} |\psi({\tt t})|d{\tt t}$$
and
$$ M_\nu(\varphi) := \max_{|{\tt h}|= \nu}M^{\nu}_{\tt h}(\varphi)\quad \widetilde{M}_{\nu}(\psi):= \max_{|{\tt h}|= \nu}\widetilde{M}^{\nu}_{\tt h}(\psi).$$
\vskip0.4cm
\noindent
We suppose that the following assumptions hold
\begin{enumerate}
\item[i)] for every ${\tt u}\in \mathbb{R}^N,$~ we have 
$$m_{\tt 0}^0 (\varphi,{\tt u})  = \sum_{{\tt k}\in Z\!\!\!\!Z^N} \varphi({\tt u}-{\tt k}) = 1,~\quad  \widetilde{m}^{0}_{\tt 0}(\psi) =\int_{\mathbb{R}^N} \psi({\tt t})d{\tt t}=1$$ 
\item[ii)]
for every ${\tt h} \in \mathbb{N}^N $~ and $\nu=1, \ldots, r$ we have, 
$$m^\nu_{{\tt h}}(\varphi, {\tt u}) =: m^\nu_{{\tt h}}(\varphi),$$
where $m^\nu_{{\tt h}}(\varphi)$ is a fixed real number independent of ${\tt u}.$ 
\item[iii)]
$M_r(\varphi) + \widetilde{M}_{\nu}(\psi)< +\infty$
and there holds
$$ \lim_{w\rightarrow +\infty} \sum_{\|{\tt u}-{\tt k}\|>w } |\varphi({\tt u}-{\tt k})| \|{\tt u}-{\tt k}\|^r =0$$
uniformly with respect to ${\tt u}\in \mathbb{R}^N.$
\end{enumerate}

\vskip0.4cm
\noindent
{\bf Remark 1}
\begin{enumerate}
\item[1.]
Note that for any $\nu\in \mathbb{N}_0$ for the absolute moment $M_\nu(\varphi)$ we obtain
$$M_\nu(\varphi)\leq \sup_{{\tt u}\in \mathbb{R}^N} \sum_{{\tt k}\in \mathbb{Z}^N} |\varphi({\tt u}-{\tt k})| \|{\tt u}-{\tt k}\|^{\tt \nu}\leq N^{\nu}M_\nu(\varphi).$$
Indeed for $\nu\in \mathbb{N}_0$ and ${\tt h}\in \mathbb{N}_0^N,~{\tt h}=(h_1, \ldots, h_N), $~with $|{\tt h}|= \nu$ for the first inequality we have
$$[{\tt u}-{\tt k}]^{\tt h}= |u_1 - k_1|^{h_1} \ldots |u_N - k_N|^{h_N} \leq \|{\tt u}-{\tt k}\|^{h_1}\ldots \|{\tt u}-{\tt k}\|^{h_N} =\|{\tt u}-{\tt k}\|^{\nu}.$$
For the second inequality, taking into account the elementary inequality
$$\|{\tt u} - {\tt k}\|^\nu \leq \bigg(\sum_{i=1}^N |u_i - k_i|\bigg)^\nu \leq N^{\nu-1}\sum_{i=1}^N |u_i - k_i|^\nu,$$
we have
\begin{eqnarray*}
&&\sup_{{\tt u}\in \mathbb{R}^N} \sum_{{\tt k}\in \mathbb{Z}^N} |\varphi({\tt u}-{\tt k})| \|{\tt u}-{\tt k}\|^{\tt \nu}\leq 
N^{\nu-1} \sup_{{\tt u}\in \mathbb{R}^N} \sum_{{\tt k}\in \mathbb{Z}^N} |\varphi({\tt u}-{\tt k})| \bigg(\sum_{i=1}^N |u_i - k_i|^\nu \bigg)\\&=&
N^{\nu-1} \sum_{i=1}^N \bigg(\sup_{{\tt u}\in \mathbb{R}^N} \sum_{{\tt k}\in \mathbb{Z}^N} |\varphi({\tt u}-{\tt k})| |u_i - k_i|^\nu\bigg)\leq 
N^{\nu-1} \sum_{i=1}^N  M_\nu(\varphi)= N^\nu M_\nu(\varphi).
\end{eqnarray*}
Analogously we have also
$$ \widetilde{M}_\nu (\psi) \leq \int_{\mathbb{R}^N} |\psi({\tt t})| \|{\tt t}\|^{\nu} d{\tt t}\leq N^\nu \widetilde{M}_\nu(\psi).$$
\item[2.]
Note that for $\mu,~\nu \in \mathbb{N}_0$~with $\mu< \nu,$~
$M_\nu(\varphi)< +\infty$~ implies $M_\mu(\varphi) <+\infty.$~
When the function $\varphi$ has compact support, we immediately have that $M_\nu(\varphi) < +\infty,$ for every $\nu \in \mathbb{N}_0.$
\item[3.]
Assumption iii) implies that for $j=0,1, \ldots, r-1$ there holds
$$ \lim_{w\rightarrow +\infty} \sum_{\|{\tt u}-{\tt k}\|>w } |\varphi({\tt u}-{\tt k})| \|{\tt u}-{\tt k}\|^j =0$$
uniformly with respect to $u\in I\!\!R^N.$
Indeed, for example 
$$  \sum_{\|{\tt u}-{\tt k}\|>w } |\varphi({\tt u}-{\tt k})| < \frac{1}{w^r} \sum_{\|{\tt u}-{\tt k}\|>w } |\varphi({\tt u}-{\tt k})| \|{\tt u}- {\tt k}\|^r.$$
\end{enumerate}

\vskip0.3cm
Below, the following Taylor formula for functions $ f \in C^r$ locally at the point 
${\tt x},$ will be used
\begin{eqnarray*}
f({\tt u}) &=& f({\tt x}) + \sum_{\nu =1}^r \sum_{|{\tt h}| =\nu} \frac{D^{{\tt h}}f({\tt x})}{{\tt h}!} ({\tt u} - {\tt x})^{\tt h}+ \lambda({\tt u} - {\tt x}) \|{\tt u} - {\tt x}\|^r,
\end{eqnarray*}
where $\lambda$ is a bounded function such that $\lim_{{\tt v} \rightarrow {\tt 0}} \lambda({\tt v}) = 0.$

Let ${\tt p}: \mathbb{N} \rightarrow \mathbb{R}^N,$ ${\tt p}= (p_1, \ldots, p_N)$ be a function such that ${\tt p}(n)>{\tt 0}$ for every $n>0$~ and satisfying the following assumptions
\begin{enumerate}
\item[j)] for every $i=1, \ldots, N,$~ $\lim_{n \rightarrow +\infty} p_i(n) = +\infty,$
\item[jj)] for every $i=1, \ldots, N$ it holds
$$ \lim_{n\rightarrow +\infty} \frac{n}{p_i(n)} =: a_i >0.$$
Putting ${\tt a} = (a_1, \ldots, a_N)$ assumption jj) can be read as 
$$\lim_{n\rightarrow +\infty} \frac{n}{{\tt p}(n)} = {\tt a}.$$
\end{enumerate}

For $n \in \mathbb{N}, \varphi$ and $\psi$ satisfying the above assumptions, the multivariate generalized Durrmeyer sampling type series generated by $\varphi$ and $\psi$ is defined as 
$$(S_{{\tt p}(n)}^{\varphi, \psi}f)({\tt x}) = \sum_{{\tt k} \in \mathbb{Z}^N}  \varphi({\tt p}(n){\tt x} - {\tt k}) \bigg[ \langle{\tt p}(n)\rangle
 \int_{\mathbb{R}^N} \psi({\tt p}(n){\tt u} - {\tt k}) f({\tt u})du \bigg]\quad {\tt x} \in \mathbb{R}^N.$$
We will put $Dom S = \bigcap_{n \in I\!\!N} Dom S_{{\tt p}(n)}^{\varphi, \psi}$~ where $Dom S_{{\tt p}(n)}^{\varphi, \psi}$~ is the space of all functions $f:\mathbb{R}^N \rightarrow \mathbb{R}$~ for which the series defining $S_{{\tt p}(n)}^{\varphi, \psi}$~ is absolutely convergent for every 
${\tt x} \in \mathbb{R}^N.$~ 
Under the above assumptions it is easy to see that  $Dom S_{{\tt p}(n)}^{\varphi, \psi}$~ contains the space $L^\infty (\mathbb{R}^N).$
In the particular case when $\varphi$~ has compact support, every function $f$~ belongs to $Dom S.$ Indeed
for every fixed $n$~and $ {\tt x}\in \mathbb{R}^N,$~ only a finite numbers of ${\tt k}$~ occurs in the series defining $(S_{{\tt p}(n)}^{\varphi, \psi}f).$~

\section{Asymptotic behaviour and Voronovskaja formula}
In order to study the asymptotic behaviour of the multivariate Durrmeyer type operator, we need the following further notations. For given vectors ${\tt h} = (h_1, \ldots, h_N), {\tt j} = (j_1, \ldots, j_N)  \in \mathbb{N}^N, {\tt a}= (a_1, \ldots a_N)\in \mathbb{R}^N,$ we set
$$\left(\begin{array}{c} {\tt h}\\ {\tt j}\end{array} \right) = \left(\begin{array}{c} h_1\\ j_1\end{array} \right)\cdots \left(\begin{array}{c} h_N\\ j_N\end{array} \right),$$
$$\sum_{{\tt j}={\tt 0}}^{\tt h} 
\left(\begin{array}{c} {\tt h}\\ {\tt j}\end{array} \right) {\tt a}^{\tt j} = 
\sum_{ j_1= 0}^{ h_1}\ldots  \sum_{j_N= 0}^{ h_N} 
\left(\begin{array}{c} h_1\\ j_1\end{array} \right)\cdots \left(\begin{array}{c} h_N\\ j_N\end{array} \right)a_1^{j_1}\cdots 
a_N^{j_N}.$$
If $A_{\tt j}$ is scalar we set
$$\sum_{{\tt j}={\tt 0}}^{\tt h} 
\left(\begin{array}{c} {\tt h}\\ {\tt j}\end{array} \right) A_{\tt j} = 
\sum_{ j_1= 0}^{ h_1}\left(\begin{array}{c} h_1\\ j_1\end{array} \right)\ldots  \sum_{j_N= 0}^{ h_N} \left(\begin{array}{c} h_N\\ j_N\end{array} \right)A_{(j_1, \cdots, j_N)}.$$
We have the following theorem.
\newtheorem{Theorem}{Theorem}
\begin{Theorem}\label{pointwise}
Let $f \in L^\infty(\mathbb{R}^N)$~ and let ${\tt x} \in \mathbb{R}^N$~ be fixed.
Assume that $f$~ belongs to $C^r$~ locally at ${\tt x}.$~
Under the assumptions i), ii), iii) we have, for $~n\rightarrow + \infty,$
\begin{eqnarray*}
 && (S_{{\tt p}(n)}^{\varphi, \psi}f)({\tt x}) - f({\tt x}) = \sum_{\nu = 1}^r 
\sum_{|{\tt h}| = \nu} \frac{D^{{\tt h}}f({\tt x})}{{\tt h}!({\tt p}(n))^{{\tt h}}} 
\sum_{{\tt s}= {\tt 0}}^{\tt h} \left(\begin{array}{c} {\tt h}\\ {\tt s}\end{array} \right) m^{|{\tt h-s}|}_{{\tt h-s}}(\varphi)  \widetilde{m}^{|{\tt s}|}_{\tt s}(\psi) + o(\|{\tt p}(n)\|^{-r}).
\end{eqnarray*}
\end{Theorem}
{\bf Proof.}~
Using the local Taylor formula of order $r$ for the function $f,$~ there exists a bounded function 
$\lambda$~ such that 
$\lim_{{\tt v}\rightarrow {\tt 0}} \lambda({\tt v}) = 0$~ and 
\begin{eqnarray*}
f({\tt u}) = f({\tt x}) + \sum_{\nu =1}^r \sum_{|{\tt h}| =\nu} \frac{D^{{\tt h}}f({\tt x})}{{\tt h}!}({\tt u} - {\tt x})^{\tt h} +
\lambda({\tt u} - {\tt x}) \|{\tt u} - {\tt x}\|^r.
\end{eqnarray*}
Thus, by i) and ii) we have,
\begin{eqnarray*}
&& (S_{{\tt p}(n)}^{\varphi, \psi}f)({\tt x}) - f({\tt x}) =  \sum_{{\tt k} \in \mathbb{Z}^N}\varphi ({\tt p}(n){\tt x} - {\tt k})  \bigg[ \langle{\tt p}(n)\rangle
 \int_{\mathbb{R}^N} \psi({\tt p}(n){\tt u} - {\tt k}) (f({\tt u})-f({\tt x}))du \bigg]\\&=&
 \sum_{{\tt k} \in \mathbb{Z}^N}\varphi ({\tt p}(n){\tt x} - {\tt k})  \bigg[ \langle{\tt p}(n)\rangle
 \int_{\mathbb{R}^N} \psi({\tt p}(n){\tt u} - {\tt k}) \sum_{\nu =1}^r \sum_{|{\tt h}| =\nu} \frac{D^{{\tt h}}f({\tt x})}{{\tt h}!}({\tt u} - {\tt x})^{\tt h}du \bigg]\\ &+& 
 \sum_{{\tt k} \in \mathbb{Z}^N}\varphi ({\tt p}(n){\tt x} - {\tt k})  \bigg[ \langle{\tt p}(n)\rangle
 \int_{\mathbb{R}^N} \psi({\tt p}(n){\tt u} - {\tt k}) \lambda({\tt u} - {\tt x}) \|{\tt u} - {\tt x}\|^r du \bigg]\\&=&
 I_1 + \ldots + I_r +J,
\end{eqnarray*}
where for $j=1, \cdots, r$
$$ I_j=\sum_{|{\tt h}| =j} \frac{D^{{\tt h}}f({\tt x})}{{\tt h}!} \sum_{{\tt k} \in \mathbb{Z}^N}\varphi ({\tt p}(n){\tt x} - {\tt k})\bigg[ \langle{\tt p}(n)\rangle
 \int_{\mathbb{R}^N} \psi({\tt p}(n){\tt u} - {\tt k})  ({\tt u} - {\tt x})^{\tt h}du \bigg]$$
and 
$$ J:= \sum_{{\tt k} \in \mathbb{Z}^N}\varphi ({\tt p}(n){\tt x} - {\tt k})\bigg[ \langle{\tt p}(n)\rangle
 \int_{\mathbb{R}^N} \psi({\tt p}(n){\tt u} - {\tt k}) \lambda({\tt u} - {\tt x}) \|{\tt u} - {\tt x}\|^r du \bigg].$$
Now we evaluate $I_j.$ Putting ${\tt p}(n){\tt u} - {\tt k}={\tt t}$ we obtain ${\tt u }= \frac{{\tt t} +{\tt k} }{{\tt p}(n)}$ and the Jacobian is $ \langle{\tt p}(n)\rangle^{-1}$ so that 
\begin{eqnarray*}
&& \langle{\tt p}(n)\rangle
 \int_{\mathbb{R}^N} \psi({\tt p}(n){\tt u} - {\tt k}) ({\tt u} - {\tt x})^{\tt h}du=
 \langle{\tt p}(n)\rangle \int_{\mathbb{R}^N} \psi({\tt t}) \bigg(\frac{{\tt t} +{\tt k} }{{\tt p}(n)}- {\tt x} \bigg)^{\tt h} \frac{d{\tt t}}{\langle{\tt p}(n)\rangle}\\ &=&
\frac{1}{({\tt p}(n))^{\tt h}} \int_{\mathbb{R}^N} \psi({\tt t}) ({\tt t} + {\tt k} - {\tt p}(n){\tt x})^{\tt h} d{\tt t}.
\end{eqnarray*}
Now we have 
\begin{eqnarray*}
&&({\tt t} + {\tt k} - {\tt p}(n){\tt x})^{\tt h} \\&=&
\sum_{s_1=0}^{h_1} \left(\begin{array}{c} h_1\\ s_1\end{array} \right) t_1^{s_1} (k_1 -p_1(n)x_1)^{h_1-s_1}\cdot \ldots \cdot
\sum_{s_N=0}^{h_N} \left(\begin{array}{c} h_N\\ s_N\end{array} \right) t_N^{s_N} (k_N -p_N(n)x_N)^{h_N-s_N}\\ &=&
\sum_{{\tt s}= {\tt 0}}^{\tt h} \left(\begin{array}{c} {\tt h}\\ {\tt s}\end{array} \right) {\tt t}^{\tt s} ({\tt k}- {\tt p}(n){\tt x})^{{\tt h}-{\tt s}}
\end{eqnarray*}
 so that
\begin{eqnarray*}
&& I_j = \sum_{|{\tt h}| =j} \frac{D^{{\tt h}}f({\tt x})}{{\tt h}!} \sum_{{\tt k} \in \mathbb{Z}^N}\varphi ({\tt p}(n){\tt x} - {\tt k}) \bigg[ \frac{1}{({\tt p}(n))^{{\tt h}}}\int_{\mathbb{R}^N} \psi({\tt t}) \sum_{{\tt s}= {\tt 0}}^{\tt h} \left(\begin{array}{c} {\tt h}\\ {\tt s}\end{array} \right) {\tt t}^{\tt s} ({\tt k}- {\tt p}(n){\tt x})^{{\tt h}- {\tt s}}d{\tt t}\bigg]\\ &=&
\sum_{|{\tt h}| =j} \frac{1}{({\tt p}(n))^{{\tt h}}} \frac{D^{{\tt h}}f({\tt x})}{{\tt h}!} \sum_{s=0}^{\tt h} \left(\begin{array}{c} {\tt h}\\ {\tt s}\end{array} \right)\sum_{{\tt k} \in \mathbb{Z}^N}\varphi ({\tt p}(n){\tt x} - {\tt k}) ({\tt k} - {\tt p}(n){\tt x})^{\tt h-s}\bigg[  \int_{\mathbb{R}^N} \psi({\tt t}){\tt t}^{\tt s}  d{\tt t} \bigg]\\ &=&
\sum_{|{\tt h}| =j} \frac{1}{({\tt p}(n))^{{\tt h}}} \frac{D^{{\tt h}}f({\tt x})}{{\tt h}!} \sum_{{\tt s}= {\tt 0}}^{\tt h} \left(\begin{array}{c} {\tt h}\\ {\tt s}\end{array} \right) m^{|{\tt h-s}|}_{{\tt h-s}}(\varphi)  \widetilde{m}^{|{\tt s}|}_{\tt s}(\psi).
\end{eqnarray*}
 
We now evaluate $J.$ 
Let $\varepsilon >0$ be fixed. There exists $\delta >0$ such that $|\lambda ({\tt v})| < \varepsilon$ for every $\|{\tt v}\| \leq \delta.$ 

Hence,
\begin{eqnarray*}
|J| &\leq& \sum_{{\tt k} \in \mathbb{Z}^N}|\varphi ({\tt p}(n){\tt x} - {\tt k})| \bigg[ \langle{\tt p}(n)\rangle
\bigg(\int_{\|{\tt u} - {\tt x}\| \leq \delta} + \int_{\|{\tt u} - {\tt x}\| > \delta} \bigg) |\psi ( {\tt p}(n){\tt u} - {\tt k})| |\lambda ({\tt u} - {\tt x})| \|{\tt u} - {\tt x}\|^rd{\tt u} \bigg] \\&=& J' + J''.
\end{eqnarray*}
For $J'$ we obtain, putting ${\tt p}(n){\tt u} - {\tt k}={\tt t},$ 

\begin{eqnarray*}
&&J' \leq 
\varepsilon \sum_{{\tt k} \in \mathbb{Z}^N}|\varphi ({\tt p}(n){\tt x} - {\tt k})| \bigg[  \int_{\|{\tt \frac{{\tt t} +{\tt k} }{{\tt p}(n)}} - {\tt x}\| \leq \delta} |\psi({\tt t})| \bigg\|{\tt {\tt \frac{{\tt t} +{\tt k} }{{\tt p}(n)}} } - {\tt x}\bigg\|^rd{\tt t} \bigg]\\ &=&
\frac{\varepsilon}{\|{\tt p}(n) \|^r}\sum_{{\tt k} \in \mathbb{Z}^N}|\varphi ({\tt p}(n){\tt x} - {\tt k})| \bigg[  \int_{\|{\tt \frac{{\tt t} +{\tt k} }{{\tt p}(n)}} - {\tt x}\| \leq \delta} |\psi({\tt t})| \|{\tt t} +{\tt k}  - {\tt p}(n){\tt x}\|^r d{\tt t} \bigg].
\end{eqnarray*}
Since
$$\|{\tt t} +{\tt k}  - {\tt p}(n){\tt x}\|^r \leq (\|\tt t\| + \|{\tt k}  - {\tt p}(n){\tt x}\|)^r \leq 2^{r-1} (\|{\tt t}\|^r  + \|{\tt k}  - {\tt p}(n){\tt x}\|^r ),$$
we have
\begin{eqnarray*}
&&J' \leq \frac{\varepsilon 2^{r-1}}{\|{\tt p}(n) \|^r}\sum_{{\tt k} \in \mathbb{Z}^N}|\varphi ({\tt p}(n){\tt x} - {\tt k})| \\&\times& \bigg[\int_{\|{\tt \frac{{\tt t} +{\tt k} }{{\tt p}(n)}} - {\tt x}\| \leq \delta} |\psi({\tt t})| \|{\tt t}\|^rd{\tt t} + \int_{\|{\tt \frac{{\tt t} +{\tt k} }{{\tt p}(n)}} - {\tt x}\| \leq \delta} |\psi({\tt t})|\|{\tt k}  - {\tt p}(n){\tt x}\|^rd{\tt t} \bigg]\\  &\leq&
\frac{\varepsilon 2^{r-1}}{\|{\tt p}(n) \|^r} \bigg(N^r \widetilde{M}_{r}(\psi) M_{ 0}(\varphi) + N^r\widetilde{M}_{ 0}(\psi) M_r(\varphi)\bigg).
\end{eqnarray*}

Finally, 
\begin{eqnarray*}
&&J'' \leq \|\lambda\|_\infty \sum_{{\tt k} \in \mathbb{Z}^N} |\varphi ({\tt p}(n) {\tt x} - {\tt k})| \bigg[\int_{\|{\tt {\tt \frac{{\tt t} +{\tt k} }{{\tt p}(n)}}} - {\tt x}\|> \delta} |\psi ( {\tt t})|\bigg\|{\tt {\tt \frac{{\tt t} +{\tt k} }{{\tt p}(n)}}}  - {\tt x}\bigg\|^r d{\tt t} \bigg]
\\ &\leq& 
 \frac{2^{r-1}\|\lambda\|_\infty}{\|{\tt p}(n) \|^r}\sum_{{\tt k} \in \mathbb{Z}^N}|\varphi ({\tt p}(n){\tt x} - {\tt k})|\bigg[ \int_{\|{\tt {\tt \frac{{\tt t} +{\tt k} }{{\tt p}(n)}}} - {\tt x}\|> \delta} |\psi ( {\tt t})|(\|{\tt t}\|^r+ \|{\tt k} - {\tt p}(n){\tt x}\|^r )d{\tt t} \bigg]\\ &\leq&
\frac{2^{r-1}\|\lambda\|_\infty}{\|{\tt p}(n) \|^r}\sum_{{\tt k} \in \mathbb{Z}^N}|\varphi ({\tt p}(n){\tt x} - {\tt k})|\int_{\|{\tt {\tt \frac{{\tt t} +{\tt k} }{{\tt p}(n)}}} - {\tt x}\|> \delta} |\psi ( {\tt t})|\|{\tt t}\|^rd{\tt t}\\ &+&
\frac{2^{r-1}\|\lambda\|_\infty}{\|{\tt p}(n) \|^r}\sum_{{\tt k} \in \mathbb{Z}^N}|\varphi ({\tt p}(n){\tt x} - {\tt k})| \|{\tt k} - {\tt p}(n){\tt x}\|^r \int_{\|{\tt {\tt \frac{{\tt t} +{\tt k} }{{\tt p}(n)}}} - {\tt x}\|> \delta} |\psi ( {\tt t})|d{\tt t} = J''_1 + J''_2.
\end{eqnarray*}
For $J''_1$ we have
\begin{eqnarray*}
&&J''_1 = \frac{2^{r-1}\|\lambda\|_\infty}{\|{\tt p}(n) \|^r}\sum_{\|{\tt p}(n){\tt x} - {\tt k}\|< \|{\tt p}(n)\| \delta/2} |\varphi ({\tt p}(n){\tt x} - {\tt k})|\int_{\|{\tt {\tt \frac{{\tt t} +{\tt k} }{{\tt p}(n)}}} - {\tt x}\|> \delta} |\psi ( {\tt t})|\|{\tt t}\|^rd{\tt t}\\ &+& 
\frac{2^{r-1}\|\lambda\|_\infty}{\|{\tt p}(n) \|^r}\sum_{\|{\tt p}(n){\tt x} - {\tt k}\|\geq \|{\tt p}(n)\| \delta/2} |\varphi ({\tt p}(n){\tt x} - {\tt k})|\int_{\|{\tt {\tt \frac{{\tt t} +{\tt k} }{{\tt p}(n)}}} - {\tt x}\|> \delta} |\psi ( {\tt t})|\|{\tt t}\|^rd{\tt t} = J''_{1,1} + J''_{1,2}.
\end{eqnarray*}
For $J''_{1,1},$ since $\|{\tt p}(n){\tt x} - {\tt k}\|< \|{\tt p}(n)\| \delta/2,$ we have 
$$\| {\tt t}\| + \|{\tt p}(n)\| \delta/2\geq \|{\tt t} \| + \|{\tt p}(n){\tt x} - {\tt k}\|\geq \|{\tt t}+ {\tt k}- {\tt p}(n){\tt x}\| \geq \|{\tt p}(n)\| \delta \Rightarrow \| {\tt t}\|\geq \|{\tt p}(n)\| \delta/2$$
then 
$$\int_{\|{\tt {\tt \frac{{\tt t} +{\tt k} }{{\tt p}(n)}}} - {\tt x}\|> \delta} |\psi ( {\tt t})|\|{\tt t}\|^rd{\tt t} \leq 
\int_{\|{\tt t}\|>\|{\tt p}(n)\| \delta/2} |\psi ( {\tt t})|\|{\tt t}\|^rd{\tt t}$$
and for the absolute continuity of the integral we have that for large $n$
$$\int_{\|{\tt {\tt \frac{{\tt t} +{\tt k} }{{\tt p}(n)}}} - {\tt x}\|> \delta} |\psi ( {\tt t})|\|{\tt t}\|^rd{\tt t}<\varepsilon$$
for every $k$ such that $\|{\tt p}(n){\tt x} - {\tt k}\|< \|{\tt p}(n)\| \delta/2.$
Hence $J''_{1,1}\leq \frac{2^{r-1} \varepsilon \|\lambda\|_\infty}{\|{\tt p}(n) \|^r} M_0(\varphi).$

For $J''_{1,2},$ we have for $n$ large, by condition iii)
$$ |J''_{1,2}| \leq \frac{2^{r-1}\|\lambda\|_\infty}{\|{\tt p}(n) \|^r} \widetilde{M}_r(\psi) N^r \sum_{\|{\tt p}(n){\tt x} - {\tt k}\|\geq \|{\tt p}(n)\| \delta/2} |\varphi ({\tt p}(n){\tt x} - {\tt k})|\leq \frac{2^{r-1}\|\lambda\|_\infty}{\|{\tt p}(n) \|^r} \widetilde{M}_r(\psi)N^r \varepsilon$$
so the assertion follows, estimating in an analogous way the term $J''_2.$ $\Box$

\vskip0.4cm
\noindent
As a consequence, we obtain the following Voronovskaja formula
\newtheorem{Corollary}{Corollary}
\begin{Corollary}\label{voronov}
Under the assumptions of Theorem 1, we have
\begin{description}
\item[i)]
For $r=1$ there holds
\begin{eqnarray*}
 && \lim_{n \rightarrow + \infty}n [(S_{{\tt p}(n)}^{\varphi, \psi}f)({\tt x}) - f({\tt x})] =
\sum_{|{\tt h}| = 1} {\tt a}^{{\tt h}} D^{{\tt h}}f({\tt x}) 
\sum_{{\tt s}= {\tt 0}}^{\tt h} \left(\begin{array}{c} {\tt h}\\ {\tt s}\end{array} \right) m^{|{\tt h-s}|}_{{\tt h-s}}(\varphi)  \widetilde{m}^{|{\tt s}|}_{\tt s}(\psi) \\
&&= \sum_{k=1}^N a_k\frac{\partial f}{\partial x_k}({\tt x})\{m^1_{{\tt e}_k}(\varphi) + \widetilde{m}^1_{{\tt e}_k}(\psi)\}.
\end{eqnarray*}

\item[ii)]
For $r\in \mathbb{N}, r>1$ if moreover for every ${\tt h} \in \mathbb{N}^N$, with $|{\tt h}| = \nu$ and $\nu=1, \ldots, r-1$ we have 
$ m^\nu_{{\tt h}}(\varphi) =\widetilde{m}^\nu_{{\tt h}}(\psi)= 0$ and $m^r_{{\tt h}}(\varphi), \widetilde{m}^r_{{\tt h}}(\psi) \neq 0$ for some ${\tt h}$ with $|{\tt h}|=r,$ then we have
\begin{eqnarray*}
&&\lim_{n \rightarrow + \infty}n^r [(S_{{\tt p}(n)}^{\varphi, \psi}f)({\tt x}) - f({\tt x})] = \sum_{|{\tt h}| = r} {\tt a}^{{\tt h}}\frac{D^{{\tt h}}f({\tt x})}{{\tt h}!} \sum_{{\tt s}= {\tt 0}}^{\tt h} \left(\begin{array}{c} {\tt h}\\ {\tt s}\end{array} \right) m^{|{\tt h-s}|}_{{\tt h-s}}(\varphi)  \widetilde{m}^{|{\tt s}|}_{\tt s}(\psi).
\end{eqnarray*}
\end{description}
\end{Corollary}
\vskip0.4cm
\noindent
{\bf Remark 2}. Note that for $r=2$ the relation of part ii) in Corollary \ref{voronov}, reads
\begin{eqnarray*}
&&\lim_{n \rightarrow + \infty}n^2 [(S_{{\tt p}(n)}^{\varphi, \psi}f)({\tt x}) - f({\tt x})] = \sum_{|{\tt h}| = 2} {\tt a}^{{\tt h}}\frac{D^{{\tt h}}f({\tt x})}{{\tt h}!} \sum_{{\tt s}= {\tt 0}}^{\tt h} \left(\begin{array}{c} {\tt h}\\ {\tt s}\end{array} \right) m^{|{\tt h-s}|}_{{\tt h-s}}(\varphi)  \widetilde{m}^{|{\tt s}|}_{\tt s}(\psi)\\
&=& \sum_{k=1}^N \frac{a_k^2}{2} \frac{\partial^2 f({\tt x})}{\partial x_k^2}(m^2_{2{\tt e}_k}(\varphi) + 2m^1_{{\tt e}_k}(\varphi) \widetilde{m}^1_{{\tt e}_k}(\psi) + \widetilde{m}^2_{2{\tt e}_k}(\psi))\\ &+&
\sum_{i,j=1, i\neq j}^N a_i a_j \frac{\partial^2 f({\tt x})}{\partial x_i \partial x_j} (m^2_{{\tt e}_i+{\tt e}_j}(\varphi) + m^1_{{\tt e}_j}(\varphi)\widetilde{m}^1_{{\tt e}_i}(\psi)  + m^1_{{\tt e}_i}(\varphi)\widetilde{m}^1_{{\tt e}_j}(\psi) +  \widetilde{m}^2_{{\tt e}_i + {\tt e}_j}(\psi) ).
\end{eqnarray*}

\section{Quantitative estimates}

In this section our aim is to determine the order of convergence in Corollary 1 using the classical Peetre K-functional, introduced by J. Peetre (see \cite{PE2}) and defined in the multivariate setting by
$$ K(\varepsilon, f, C^0, C^1)\equiv K(\varepsilon, f) := \inf \{ \|f-g\|_\infty + \varepsilon \max_{i=1,\cdots,N}\bigg \|\frac{\partial g}{\partial x_i} \bigg\|_\infty : g\in C^1 \}$$
for $f\in C^0$~ and $\varepsilon \geq 0,$~ (see also \cite{BM1}).

If $f$~belongs to $C^r$~we have the Taylor formula with the remainder in the Lagrange form
\begin{eqnarray*}
f({\tt u}) &=& f({\tt x}) + \sum_{\nu =1}^r \sum_{|{\tt h}| =\nu} \frac{D^{{\tt h}}f({\tt x})}{{\tt h}!} ({\tt u} - {\tt x})^{\tt h}+ R^L_r(f; {\tt x}, {\tt u}),
\end{eqnarray*}
for ${\tt x}, {\tt u}\in \mathbb{R}^N,~ r\geq 1,$~where 
$$ R^L_r(f,{\tt x},{\tt u}) = \sum_{|{\tt h}|=r} ({\tt u}-{\tt x})^{{\tt h}} \frac{D^{\tt h} f({\tt \xi})}{{\tt h}!}=\frac{1}{r!}\sum_{i_1=1}^N  \cdots \sum_{i_r=1}^N (u_{i_1}- x_{i_1})\cdots (u_{i_r}- x_{i_r})\frac{\partial^r f({\tt \xi})}{\partial x_{i_1}\partial x_{i_2}\cdots \partial x_{i_r}}.$$

Now we give an extension of Lemma 4.2 in \cite{BM1} (see also \cite{GPR} for the one-dimensional case) obtaining the following estimate of the remainder $R_r(f; {\tt x}, {\tt u})$~ 

\newtheorem{Lemma}{Lemma}
\begin{Lemma}\label{remainder}
Let $f\in C^r$~ and ${\tt x}, {\tt u} \in \mathbb{R}^N.$~
Then 
$$|R_r(f; {\tt x}, {\tt u})| \leq \frac{2}{r!} \sum_{i_1=1}^N \sum_{i_2=1}^N \cdots \sum_{i_r=1}^N|u_{i_1}- x_{i_1}|^r K \bigg(\frac{N}{2(r+1)} |u_{i_1}- x_{i_1}|, \frac{\partial^{r}f }{\partial x_{i_1}\cdots \partial x_{i_r}}\bigg).$$

\end{Lemma}
{\bf Proof.}~
We have, using the Taylor formula with the remainder in the Peano form,
$$R_r(f; {\tt x}, {\tt u}) = f({\tt u}) - f({\tt x}) - \sum_{\nu =1}^r \sum_{|{\tt h}| =\nu} \frac{D^{{\tt h}}f({\tt x})}{{\tt h}!} ({\tt u} - {\tt x})^{\tt h}.$$
Using the Taylor formula with remainder in the Lagrange form, there exists a point ${\tt \xi}\in L({\tt x}, {\tt u}),$ being $L({\tt x}, {\tt u})$ the segment with end points ${\tt x}, {\tt u},$ such that 
$$ f({\tt u}) - f({\tt x})=  \sum_{\nu =1}^{r-1} \sum_{|{\tt h}| =\nu} \frac{D^{{\tt h}}f({\tt x})}{{\tt h}!} ({\tt u} - {\tt x})^{\tt h}+  \sum_{|h|=r}( {\tt u} - {\tt x})^{\tt h}\frac{D^{{\tt h}}f({\tt \xi})}{{\tt h}!} .$$
So we have 
\begin{eqnarray*}
&&|R_r(f; {\tt x}, {\tt u})| = 
\bigg|\sum_{|h|=r}( {\tt u} - {\tt x})^{\tt h} \frac{D^{{\tt h}}f({\tt \xi})}{{\tt h}!} - \sum_{|h|=r}( {\tt u} - {\tt x})^{\tt h} \frac{D^{{\tt h}}f({\tt x})}{{\tt h}!} \bigg| \\&=&
\bigg| \sum_{|h|=r} \frac{({\tt u} - {\tt x})^{\tt h}}{{\tt h}!}  (D^{{\tt h}} f({\tt \xi}) - D^{{\tt h}}f({\tt x})) \bigg|
\\ &\leq& 
\frac{2}{r!} \sum_{i_1=1}^N \sum_{i_2=1}^N \cdots \sum_{i_r=1}^N|u_{i_1}- x_{i_1}||u_{i_2}- x_{i_2}| \cdots |u_{i_r}- x_{i_r}|\bigg \| \frac{\partial^{r}f }{\partial x_{i_1}\partial x_{i_2}\cdots \partial x_{i_r}} \bigg\|_\infty.
\end{eqnarray*}

Let now $g\in C^{r+1}$~ be fixed and
$$R^L_r(g; {\tt x}, {\tt u})= \frac{1}{(r+1)!}\sum_{i_1=1}^N  \cdots \sum_{i_{r+1}=1}^N (u_{i_1}- x_{i_1})\cdots (u_{i_{r+1}}- x_{i_{r+1}})\frac{\partial^{r+1} g({\tt \xi'})}{\partial x_{i_1}\partial x_{i_2}\cdots \partial x_{i_{r+1}}} $$
where ${\tt \xi'}\in L({\tt x}, {\tt u}).$
We get
\begin{eqnarray*}
&& |R_r(f; {\tt x}, {\tt u})| \leq
|R_r(f-g; {\tt x}, {\tt u})| + | R^L_r(g; {\tt x}, {\tt u})|\\ &\leq&
\frac{2}{r!} \sum_{i_1=1}^N \sum_{i_2=1}^N \cdots \sum_{i_r=1}^N|u_{i_1}- x_{i_1}||u_{i_2}- x_{i_2}| \cdots |u_{i_r}- x_{i_r}|\bigg \| \frac{\partial^{r}(f-g) }{\partial x_{i_1}\partial x_{i_2}\cdots \partial x_{i_r}} \bigg\|_\infty\\ &+& 
\frac{1}{(r+1)!}\sum_{i_1=1}^N  \cdots \sum_{i_{r+1}=1}^N |u_{i_1}- x_{i_1}|\cdots |u_{i_{r+1}}- x_{i_{r+1}}| \bigg\|\frac{\partial^{r+1} g}{\partial x_{i_1}\partial x_{i_2}\cdots \partial x_{i_{r+1}}}\bigg\|.
\end{eqnarray*}
Using the inequality $ n a_1 a_2 \cdots a_n \leq a_1^n + a_2^n + \cdots + a_n^n$ for $a_i>0,$ we have
\begin{eqnarray*}
&& |R_r(f; {\tt x}, {\tt u})|\leq
 \frac{2}{r!} \sum_{i_1=1}^N \sum_{i_2=1}^N \cdots \sum_{i_r=1}^N |u_{i_1}- x_{i_1}|^r \bigg \| \frac{\partial^{r}(f-g) }{\partial x_{i_1}\partial x_{i_2}\cdots \partial x_{i_r}} \bigg\|_\infty\\ &+& 
\frac{1}{(r+1)!}\sum_{i_1=1}^N \sum_{i_2=1}^N \cdots \sum_{i_{r+1}=1}^N |u_{i_1}- x_{i_1}|^{r+1} \bigg\| \frac{\partial^{r+1} g}{\partial x_{i_1}\partial x_{i_2}\cdots \partial x_{i_{r+1}}} \bigg\|_\infty
\\ &\leq&
\frac{2}{r!} \sum_{i_1=1}^N \sum_{i_2=1}^N \cdots \sum_{i_r=1}^N|u_{i_1}- x_{i_1}|^r  \\ &\times& \bigg( \bigg\| \frac{\partial^{r}(f-g) }{\partial x_{i_1}\cdots \partial x_{i_{r}}} \bigg\|_\infty + \frac{N}{2(r+1)} |u_{i_1}- x_{i_1}| \max_{i_{r+1}= 1, \cdots, N} \bigg\| \frac{\partial^{r+1} g}{\partial x_{i_1}\cdots \partial x_{i_{r+1}}} \bigg\|_\infty \bigg).
\end{eqnarray*}
Taking now the infimum over $g\in C^{r+1}$ we get
\begin{eqnarray*}
&& |R_r(f; {\tt x}, {\tt u})|\leq
\frac{2}{r!} \sum_{i_1=1}^N \sum_{i_2=1}^N \cdots \sum_{i_r=1}^N|u_{i_1}- x_{i_1}|^r K \bigg(\frac{N}{2(r+1)} |u_{i_1}- x_{i_1}|, \frac{\partial^{r}f }{\partial x_{i_1}\cdots \partial x_{i_r}}\bigg). \Box
\end{eqnarray*}
\vskip0,4cm
\noindent
In particular if $f\in C^1$~ we have (see \cite{BM2})
$$|R_1(f; {\tt x}, {\tt u})|\leq 2 \sum_{i =1}^N |x_i - u_i| K\bigg( \frac{N}{4}|x_i - u_i|, \frac{\partial f}{\partial x_i}\bigg).$$

Here we study an estimate of the convergence in Corollary \ref{voronov}.
\begin{Theorem}\label{kfunctional}
Under the assumptions of Theorem \ref{pointwise} we have 
\begin{description}
\item[i)] For $r=1,$ let $f \in C^1$~and ${\tt x}\in \mathbb{R}^N$~be fixed.  If moreover $M_2(\varphi) + \widetilde{M}_2(\psi) < +\infty$ 
then there holds 
\begin{eqnarray*}
&&\bigg|n((S_{{\tt p}(n)}^{\varphi, \psi}f)({\tt x}) - f({\tt x})) - \sum_{|{\tt h}| = 1} {\tt a}^{{\tt h}} D^{{\tt h}}f({\tt x}) 
\sum_{{\tt s}= {\tt 0}}^{\tt h} \left(\begin{array}{c} {\tt h}\\ {\tt s}\end{array} \right) m^{|{\tt h-s}|}_{{\tt h-s}}(\varphi)  \widetilde{m}^{|{\tt s}|}_{\tt s}(\psi)\bigg|\\ &\leq& \sum_{|{\tt h}| = 1}|D^{{\tt h}}f({\tt x})|\left| \bigg(\frac{n}{({\tt p}(n))^{{\tt h}}} -{\tt a}^{{\tt h}}\bigg)\right| \bigg|\sum_{{\tt s}= {\tt 0}}^{\tt h} \left(\begin{array}{c} {\tt h}\\ {\tt s}\end{array} \right) m^{|{\tt h-s}|}_{{\tt h-s}}(\varphi)  \widetilde{m}^{|{\tt s}|}_{\tt s}(\psi)\bigg|\\ &+&
\frac{2Nn}{\|{\tt p}(n)\|} A_1 \sum_{i =1}^N K\bigg(\frac{B_1}{A_1}\frac{N^2}{4\|{\tt p}(n)\|}, \frac{\partial f}{\partial x_i}\bigg),
\end{eqnarray*}
where
\begin{eqnarray*}
&& A_1:= M_0(\varphi) \widetilde{M}_1(\psi) + M_1(\varphi) \widetilde{M}_0(\psi)\\
&&B_1:= \widetilde{M}_2(\psi) M_0(\varphi)+ \widetilde{M}_0(\psi)M_2(\varphi) + 2\widetilde{M}_1(\psi)M_1(\varphi)
\end{eqnarray*}
\item[ii)] For $r\in \mathbb{N}, r>1,$
let $f \in C^r$~ and let ${\tt x}\in \mathbb{R}^N.$~
If moreover for every ${\tt h} \in \mathbb{N}^N,~|{\tt h}|=\nu,$ $\nu=1, \ldots, r-1$ we have $ m^\nu_{{\tt h}}(\varphi) =\widetilde{m}^\nu_{{\tt h}}(\psi)= 0$ and $m^r_{{\tt h}}(\varphi), \widetilde{m}^r_{{\tt h}}(\psi) \neq 0$ for some ${\tt h}~|{\tt h}|=r,$ and if moreover $M_{r+1}(\varphi) + \widetilde{M}_{r+1}(\psi) < +\infty,$ then there holds
\begin{eqnarray*}
&&\bigg|n^r((S_{{\tt p}(n)}^{\varphi, \psi}f)({\tt x}) - f({\tt x})) - \sum_{|{\tt h}| = r} {\tt a}^{{\tt h}}\frac{D^{{\tt h}}f({\tt x})}{{\tt h}!} \sum_{{\tt s}= {\tt 0}}^{\tt h} \left(\begin{array}{c} {\tt h}\\ {\tt s}\end{array} \right) m^{|{\tt h-s}|}_{{\tt h-s}}(\varphi)  \widetilde{m}^{|{\tt s}|}_{\tt s}(\psi)\bigg|\\ &\leq& 
\sum_{|{\tt h}| = r} \bigg |\frac{D^{{\tt h}}f({\tt x})}{\tt h!}\bigg| \left| \bigg(\frac{n^r}{({\tt p}(n))^{{\tt h}}} -{\tt a}^{\tt h}\bigg)\right| \bigg|\sum_{{\tt s}= {\tt 0}}^{\tt h} \left(\begin{array}{c} {\tt h}\\ {\tt s}\end{array} \right) m^{|{\tt h-s}|}_{{\tt h-s}}(\varphi)  \widetilde{m}^{|{\tt s}|}_{\tt s}(\psi)\bigg|\\ &+&
\frac{2}{r!}\frac{n^r N^rA_r}{\|{\tt p}(n)\|^r} \sum_{i_1=1}^N \sum_{i_2=1}^N \cdots \sum_{i_r=1}^N 
 K \bigg( \frac{N^{2}}{2(r+1) \|{\tt p}(n)\|} \frac{B_r}{A_r}, \frac{\partial^{r}f }{\partial x_{i_1}\cdots \partial x_{i_{r}}} \bigg),
\end{eqnarray*}
where 
\begin{eqnarray*}
A_r:= \sum_{\mu=0}^r \left(\begin{array}{c} r\\ \mu\end{array} \right) M_\mu(\varphi) \widetilde{M}_{r-\mu}(\psi),~~
B_r:= \sum_{\mu=0}^{r+1} \left(\begin{array}{c} r+1\\ \mu\end{array} \right)  M_\mu(\varphi) \widetilde{M}_{r+1-\mu}(\psi).
\end{eqnarray*}
\end{description}
\end{Theorem}
{\bf Proof.}~ For part i), using the Taylor formula of the first order, we have
$$f({\tt u}) = f({\tt x}) + \sum_{|{\tt h}| =1} D^{{\tt h}}f({\tt x}) ({\tt u} - {\tt x})^{\tt h}+ \lambda({\tt u} - {\tt x})\|{\tt u} - {\tt x}\|,$$
where $\lambda$ is a bounded function and $\lambda({\tt v}) \rightarrow 0$ for ${\tt v} \rightarrow {\tt 0}.$  
Following the proof of Theorem \ref{pointwise}, we obtain
\begin{eqnarray*}
&&(S_{{\tt p}(n)}^{\varphi, \psi}f)({\tt x}) - f({\tt x}) 
 =  \sum_{|{\tt h}| =1} \frac{D^{{\tt h}}f({\tt x})}{({\tt p}(n))^{{\tt h}}}  \sum_{{\tt s}= {\tt 0}}^{\tt h} \left(\begin{array}{c} {\tt h}\\ {\tt s}\end{array} \right) m^{|{\tt h-s}|}_{{\tt h-s}}(\varphi)  \widetilde{m}^{|{\tt s}|}_{\tt s}(\psi)\\&+&
\sum_{{\tt k} \in \mathbb{Z}^N}\varphi ({\tt p}(n){\tt x} - {\tt k})  \bigg[ \langle{\tt p}(n)\rangle
 \int_{\mathbb{R}^N} \psi({\tt p}(n){\tt u} - {\tt k}) \lambda({\tt u} - {\tt x}) \|{\tt u} - {\tt x}\| d{\tt u} \bigg] .
\end{eqnarray*}

Then for $f\in C^1$ we have
\begin{eqnarray*}
&&\left|(S_{{\tt p}(n)}^{\varphi, \psi}f)({\tt x}) - f({\tt x}) - \sum_{|{\tt h}| = 1} \frac{{\tt a}^{{\tt h}}}{n}D^{{\tt h}}f({\tt x}) 
\sum_{{\tt s}= {\tt 0}}^{\tt h} \left(\begin{array}{c} {\tt h}\\ {\tt s}\end{array} \right) m^{|{\tt h-s}|}_{{\tt h-s}}(\varphi)  \widetilde{m}^{|{\tt s}|}_{\tt s}(\psi)\right|
\\ &\leq& 
\left| \sum_{|{\tt h}| = 1}D^{{\tt h}}f({\tt x}) \bigg(\frac{1}{({\tt p}(n))^{{\tt h}}} -\frac{{\tt a}^{\tt h}}{n}\bigg) \sum_{{\tt s}= {\tt 0}}^{\tt h} \left(\begin{array}{c} {\tt h}\\ {\tt s}\end{array} \right) m^{|{\tt h-s}|}_{{\tt h-s}}(\varphi)  \widetilde{m}^{|{\tt s}|}_{\tt s}(\psi)\right|\\ &+&
\sum_{{\tt k} \in \mathbb{Z}^N} |\varphi ({\tt p}(n){\tt x} - {\tt k})|  \bigg[ \langle{\tt p}(n)\rangle
 \int_{\mathbb{R}^N} |\psi({\tt p}(n){\tt u} - {\tt k})| |\lambda({\tt u} - {\tt x})| \|{\tt u} - {\tt x}\| d{\tt u} \bigg]= J_1+ J_2.
\end{eqnarray*}
For $J_1$ we have immediately 
$$J_1 \leq \frac{1}{n} \sum_{|{\tt h}| = 1}|D^{{\tt h}}f({\tt x})| \left| \bigg(\frac{n}{({\tt p}(n))^{{\tt h}}} -{\tt a}^{\tt h}\bigg)\right| \bigg|\sum_{{\tt s}= {\tt 0}}^{\tt h} \left(\begin{array}{c} {\tt h}\\ {\tt s}\end{array} \right) m^{|{\tt h-s}|}_{{\tt h-s}}(\varphi)  \widetilde{m}^{|{\tt s}|}_{\tt s}(\psi)\bigg|.$$
Now we consider $J_2.$
Putting $R_1(f; {\tt x}, {\tt u}) = \lambda({\tt u} - {\tt x}) \|{\tt u} - {\tt x}\|$ by Lemma \ref{remainder} with $r=1$
we have 
$$|R_1(f; {\tt x}, {\tt u})|\leq 2 \sum_{i =1}^N |x_i - u_i| K\bigg( \frac{N}{4}|x_i - u_i|, \frac{\partial f}{\partial x_i}\bigg)  $$
and hence
$$J_2 \leq
2\sum_{{\tt k} \in \mathbb{Z}^N} |\varphi ({\tt p}(n){\tt x} - {\tt k})|  \bigg[ \langle{\tt p}(n)\rangle
 \int_{\mathbb{R}^N} |\psi({\tt p}(n){\tt u} - {\tt k})| \sum_{i =1}^N |x_i - u_i| K\bigg( \frac{N}{4}|x_i - u_i|, \frac{\partial f}{\partial x_i}\bigg)  d{\tt u} \bigg].$$
Let now $g\in C^2.$~ We have 
\begin{eqnarray*}
&& J_2 \leq 
 2\sum_{i =1}^N\bigg\|\frac{\partial (f-g)}{\partial x_i}\bigg\|_\infty \sum_{{\tt k} \in \mathbb{Z}^N} |\varphi ({\tt p}(n){\tt x} - {\tt k})|\bigg[ \langle{\tt p}(n)\rangle
 \int_{\mathbb{R}^N} |\psi({\tt p}(n){\tt u} - {\tt k})|  |x_i - u_i|d{\tt u} \bigg]\\ &+&
 \frac{N}{2}\sum_{i =1}^N  \max_{j=1,\cdots,N}\bigg \|\frac{\partial^2 g}{\partial x_i\partial x_j} \bigg\|_\infty \sum_{{\tt k} \in \mathbb{Z}^N} |\varphi ({\tt p}(n){\tt x} - {\tt k})|\bigg[ \langle{\tt p}(n)\rangle
 \int_{\mathbb{R}^N} |\psi({\tt p}(n){\tt u} - {\tt k})|  |x_i - u_i|^2  d{\tt u} \bigg]\\ &=& J_2^1 + J_2^2.
\end{eqnarray*}
For $J_2^1$ since 
$$ |x_i - u_i|\leq \|{\tt u} - {\tt x}\| \leq \frac{1}{\|{\tt p}(n)\|} \|{\tt p}(n){\tt u}- {\tt k}\| + \frac{1}{\|{\tt p}(n)\|} \|{\tt k}- {\tt p}(n) {\tt x}\| $$
using the change of variable 
 ${\tt p}(n){\tt u} - {\tt k}= {\tt t},$ we obtain
\begin{eqnarray*}
&& J_2^1\leq \frac{2}{\|{\tt p}(n)\|}\sum_{i =1}^N\bigg\|\frac{\partial (f-g)}{\partial x_i}\bigg\|_\infty \sum_{{\tt k} \in \mathbb{Z}^N} |\varphi ({\tt p}(n){\tt x} - {\tt k})| \int_{\mathbb{R}^N} |\psi({\tt t})| \|{\tt t}\| d{\tt t}\\ &+&
 \frac{2}{\|{\tt p}(n)\|}\sum_{i =1}^N\bigg\|\frac{\partial (f-g)}{\partial x_i}\bigg\|_\infty \sum_{{\tt k} \in \mathbb{Z}^N} |\varphi ({\tt p}(n){\tt x} - {\tt k})|\|{\tt k}- {\tt p}(n) {\tt x}\| \int_{\mathbb{R}^N} |\psi({\tt t})| d{\tt t}\\ &\leq&
\frac{2N}{\|{\tt p}(n)\|}\sum_{i =1}^N\bigg\|\frac{\partial (f-g)}{\partial x_i}\bigg\|_\infty ( M_0(\varphi)\widetilde{M}_1(\psi) + M_1(\varphi)\widetilde{M}_0(\psi)). 
\end{eqnarray*}
For $J_2^2$ since 
\begin{eqnarray*}
&& |x_i - u_i|^2\leq \|{\tt u} - {\tt x}\|^2 \\&\leq& \frac{1}{\|{\tt p}(n)\|^2} \|{\tt p}(n){\tt u}- {\tt k}\|^2 + \frac{1}{\|{\tt p}(n)\|^2} \|{\tt k}- {\tt p}(n) {\tt x}\|^2+ \frac{2}{\|{\tt p}(n)\|^2} \|{\tt p}(n){\tt u}- {\tt k}\|\|{\tt k}- {\tt p}(n) {\tt x}\|,
\end{eqnarray*}
we have analogously
\begin{eqnarray*}
&& J_2^2 \leq
\frac{N}{2\|{\tt p}(n)\|^2}\sum_{i =1}^N  \max_{j=1,\cdots,N}\bigg \|\frac{\partial^2 g}{\partial x_i\partial x_j} \bigg\|_\infty \sum_{{\tt k} \in \mathbb{Z}^N} |\varphi ({\tt p}(n){\tt x} - {\tt k})| \bigg[
 \int_{\mathbb{R}^N} |\psi({\tt t})| \|{\tt t}\|^2   d{\tt t} \bigg]\\ &+&
 \frac{N}{2\|{\tt p}(n)\|^2}\sum_{i =1}^N  \max_{j=1,\cdots,N}\bigg \|\frac{\partial^2 g}{\partial x_i\partial x_j} \bigg\|_\infty \sum_{{\tt k} \in \mathbb{Z}^N} |\varphi ({\tt p}(n){\tt x} - {\tt k})| \|{\tt k}- {\tt p}(n) {\tt x}\|^2 \bigg[
 \int_{\mathbb{R}^N} |\psi({\tt t})|   d{\tt t} \bigg]\\ &+&
 \frac{N}{\|{\tt p}(n)\|^2}\sum_{i =1}^N  \max_{j=1,\cdots,N}\bigg \|\frac{\partial^2 g}{\partial x_i\partial x_j} \bigg\|_\infty \sum_{{\tt k} \in \mathbb{Z}^N} |\varphi ({\tt p}(n){\tt x} - {\tt k})| \|{\tt k}- {\tt p}(n) {\tt x}\| \bigg[ 
 \int_{\mathbb{R}^N} |\psi({\tt t})| \|{\tt t}\|  d{\tt t} \bigg]\\ &\leq&
\frac{N^3}{2\|{\tt p}(n)\|^2}\sum_{i =1}^N  \max_{j=1,\cdots,N}\bigg \|\frac{\partial^2 g}{\partial x_i\partial x_j} \bigg\|_\infty (M_0(\varphi)\widetilde{M}_2(\psi) + M_2(\varphi)\widetilde{M}_0(\psi) + 2M_1(\varphi)\widetilde{M}_1(\psi)) .
\end{eqnarray*}
Then we have
\begin{eqnarray*}
&& J_2\leq \frac{2N}{\|{\tt p}(n)\|}\sum_{i =1}^N\bigg\|\frac{\partial (f-g)}{\partial x_i}\bigg\|_\infty ( M_0(\varphi)\widetilde{M}_1(\psi) + M_1(\varphi)\widetilde{M}_0(\psi))\\ &+&
\frac{N^3}{2\|{\tt p}(n)\|^2}\sum_{i =1}^N  \max_{j=1,\cdots,N}\bigg \|\frac{\partial^2 g}{\partial x_i\partial x_j} \bigg\|_\infty (M_0(\varphi)\widetilde{M}_2(\psi) + M_2(\varphi)\widetilde{M}_0(\psi) + 2M_1(\varphi)\widetilde{M}_1(\psi))\\ &=&
\frac{2N}{\|{\tt p}(n)\|} A_1 \sum_{i =1}^N \bigg( \bigg\|\frac{\partial (f-g)}{\partial x_i}\bigg\|_\infty + \frac{B_1}{A_1}\frac{N^2}{4\|{\tt p}(n)\|}  \max_{j=1,\cdots,N}\bigg \|\frac{\partial^2 g}{\partial x_i\partial x_j} \bigg\|_\infty\bigg)
\end{eqnarray*}
where $A_1= M_0(\varphi)\widetilde{M}_1(\psi) + M(\varphi)_1\widetilde{M}_0(\psi)$ and $B_1= M_0(\varphi)\widetilde{M}_2(\psi) + M_2(\varphi)\widetilde{M}_0(\psi) + 2M_1(\varphi)\widetilde{M}_1(\psi).$

Taking the infimum over $g\in C^2$~ we have
\begin{eqnarray*}
&& J_2\leq \frac{2N}{\|{\tt p}(n)\|} A_1 \sum_{i =1}^N K\bigg(\frac{B_1}{A_1}\frac{N^2}{4\|{\tt p}(n)\|}, \frac{\partial f}{\partial x_i}\bigg).
\end{eqnarray*}

For part ii), using the local Taylor formula of order r, by Theorem \ref{pointwise} and by the hypotheses
\begin{eqnarray*}
&& (S_{{\tt p}(n)}^{\varphi, \psi}f)({\tt x}) - f({\tt x})= \sum_{|{\tt h}| =r} \frac{1}{({\tt p}(n))^{{\tt h}}} \frac{D^{{\tt h}}f({\tt x})}{{\tt h}!} \sum_{{\tt s}= {\tt 0}}^{\tt h} \left(\begin{array}{c} {\tt h}\\ {\tt s}\end{array} \right) m^{|{\tt h-s}|}_{{\tt h-s}}(\varphi)  \widetilde{m}^{|{\tt s}|}_{\tt s}(\psi)\\ &+&
\sum_{{\tt k} \in \mathbb{Z}^N}\varphi ({\tt p}(n){\tt x} - {\tt k})\bigg[ \langle{\tt p}(n)\rangle
 \int_{\mathbb{R}^N} \psi({\tt p}(n){\tt u} - {\tt k}) \lambda({\tt u} - {\tt x}) \|{\tt u} - {\tt x}\|^r du \bigg].
\end{eqnarray*}
Then for $f\in C^r$ we have
\begin{eqnarray*}
&&\left|(S_{{\tt p}(n)}^{\varphi, \psi}f)({\tt x}) - f({\tt x}) - \sum_{|{\tt h}| =r} \frac{{\tt a}^{\tt h}}{n^r} \frac{D^{{\tt h}}f({\tt x})}{\tt h !} 
\sum_{{\tt s}= {\tt 0}}^{\tt h} \left(\begin{array}{c} {\tt h}\\ {\tt s}\end{array} \right) m^{|{\tt h-s}|}_{{\tt h-s}}(\varphi)  \widetilde{m}^{|{\tt s}|}_{\tt s}(\psi)\right|
\\ &\leq& 
\left| \sum_{|{\tt h}| = r}\frac{D^{{\tt h}}f({\tt x})}{\tt h !} \bigg(\frac{1}{({\tt p}(n))^{{\tt h}}} -\frac{{\tt a}^{\tt h}}{n^r}\bigg) \sum_{{\tt s}= {\tt 0}}^{\tt h} \left(\begin{array}{c} {\tt h}\\ {\tt s}\end{array} \right) m^{|{\tt h-s}|}_{{\tt h-s}}(\varphi)  \widetilde{m}^{|{\tt s}|}_{\tt s}(\psi)\right|\\ &+&
\sum_{{\tt k} \in \mathbb{Z}^N} |\varphi ({\tt p}(n){\tt x} - {\tt k})|  \bigg[ \langle{\tt p}(n)\rangle
 \int_{\mathbb{R}^N} |\psi({\tt p}(n){\tt u} - {\tt k})| |\lambda({\tt u} - {\tt x})| \|{\tt u} - {\tt x}\|^r d{\tt u} \bigg]= J_1+ J_2.
\end{eqnarray*}
For $J_1$ we have immediately 
$$J_1 \leq \frac{1}{n^r} \sum_{|{\tt h}| = r} \bigg |\frac{D^{{\tt h}}f({\tt x})}{\tt h!}\bigg| \left| \bigg(\frac{n^r}{({\tt p}(n))^{{\tt h}}} -{\tt a}^{\tt h}\bigg)\right| \bigg|\sum_{{\tt s}= {\tt 0}}^{\tt h} \left(\begin{array}{c} {\tt h}\\ {\tt s}\end{array} \right) m^{|{\tt h-s}|}_{{\tt h-s}}(\varphi)  \widetilde{m}^{|{\tt s}|}_{\tt s}(\psi)\bigg|.$$
Now we consider $J_2.$

Putting $R_r(f; {\tt x}, {\tt u}) = \lambda({\tt u} - {\tt x}) \|{\tt u} - {\tt x}\|^r$ by Lemma \ref{remainder} we have 
$$|R_r(f; {\tt x}, {\tt u})|\leq \frac{2}{r!} \sum_{i_1=1}^N \sum_{i_2=1}^N \cdots \sum_{i_r=1}^N|u_{i_1}- x_{i_1}|^r K \bigg(\frac{N}{2(r+1)} |u_{i_1}- x_{i_1}|, \frac{\partial^{r}f }{\partial x_{i_1}\cdots \partial x_{i_r}}\bigg)$$
and hence
\begin{eqnarray*}
&&J_2 \leq
\frac{2}{r!} \sum_{{\tt k} \in \mathbb{Z}^N} |\varphi ({\tt p}(n){\tt x} - {\tt k})| 
\bigg[ \langle{\tt p}(n)\rangle \int_{\mathbb{R}^N} |\psi({\tt p}(n){\tt u} - {\tt k})| \\&\times&
\sum_{i_1=1}^N \sum_{i_2=1}^N \cdots \sum_{i_r=1}^N|u_{i_1}- x_{i_1}|^r K \bigg(\frac{N}{2(r+1)} |u_{i_1}- x_{i_1}|, \frac{\partial^{r}f }{\partial x_{i_1}\cdots \partial x_{i_r}}\bigg)  d{\tt u} \bigg].
\end{eqnarray*}
Let now $g\in C^{r+1}.$~ By analogous reasonings as in part i), we have 
\begin{eqnarray*}
&&J_2 \leq
\frac{2}{r!}\sum_{i_1=1}^N \sum_{i_2=1}^N \cdots \sum_{i_r=1}^N \bigg\| \frac{\partial^{r}(f-g) }{\partial x_{i_1}\cdots \partial x_{i_{r}}} \bigg\|_\infty \\ &\times& \sum_{{\tt k} \in \mathbb{Z}^N} |\varphi ({\tt p}(n){\tt x} - {\tt k})| \bigg[ \langle{\tt p}(n)\rangle \int_{\mathbb{R}^N} |\psi({\tt p}(n){\tt u} - {\tt k})| |u_{i_1}- x_{i_1}|^r d{\tt u} \bigg]\\ &+&
\frac{N}{(r+1)!}\sum_{i_1=1}^N \sum_{i_2=1}^N \cdots \sum_{i_r=1}^N \max_{i_{r+1}= 1, \cdots, N} \bigg\| \frac{\partial^{r+1} g}{\partial x_{i_1}\cdots \partial x_{i_{r+1}}} \bigg\|_\infty\\ &\times&  \sum_{{\tt k} \in \mathbb{Z}^N} |\varphi ({\tt p}(n){\tt x} - {\tt k})|\bigg[ \langle{\tt p}(n)\rangle \int_{\mathbb{R}^N} |\psi({\tt p}(n){\tt u} - {\tt k})| |u_{i_1}- x_{i_1}|^{r+1} d{\tt u} \bigg]= J_2^1+ J_2^2.
\end{eqnarray*}
 As to $J_2^1,$ since
\begin{eqnarray*}
&& |u_{i_1}- x_{i_1}|^r \leq \|{\tt u} - {\tt x}\|^r \leq \frac{1}{\|{\tt p}(n)\|^r} (\|{\tt p}(n){\tt u}- {\tt k}\| +  \|{\tt k}- {\tt p}(n) {\tt x}\|)^r\\ &=& \frac{1}{\|{\tt p}(n)\|^r} \sum_{\mu=0}^r \left(\begin{array}{c} r\\ \mu\end{array} \right)\|{\tt p}(n){\tt u}- {\tt k}\|^{r-\mu}\|{\tt k}- {\tt p}(n) {\tt x}\|^\mu
\end{eqnarray*}
using again the change of variables ${\tt p}(n){\tt u} - {\tt k}= {\tt t},$ we obtain
\begin{eqnarray*}
&& J_2^1 \leq \frac{2}{r!}\frac{1}{\|{\tt p}(n)\|^r}  \sum_{i_1=1}^N \sum_{i_2=1}^N \cdots \sum_{i_r=1}^N \bigg\| \frac{\partial^{r}(f-g) }{\partial x_{i_1}\cdots \partial x_{i_{r}}} \bigg\|_\infty\\ &\times& \sum_{\mu=0}^r \left(\begin{array}{c} r\\ \mu\end{array} \right) \sum_{{\tt k} \in \mathbb{Z}^N} |\varphi ({\tt p}(n){\tt x} - {\tt k})| \|{\tt k}- {\tt p}(n) {\tt x}\|^\mu 
\int_{\mathbb{R}^N} |\psi({\tt t})| \|{\tt t} \|^{r-\mu} d{\tt t}\\ &\leq&
\frac{2}{r!}\frac{1}{\|{\tt p}(n)\|^r} \sum_{i_1=1}^N \sum_{i_2=1}^N \cdots \sum_{i_r=1}^N \bigg\| \frac{\partial^{r}(f-g) }{\partial x_{i_1}\cdots \partial x_{i_{r}}} \bigg\|_\infty \sum_{\mu=0}^r \left(\begin{array}{c} r\\ \mu\end{array} \right)
N^\mu M_\mu(\varphi) N^{r-\mu}\widetilde{M}_{r-\mu}(\psi)\\ &=&
\frac{2}{r!}\frac{N^r}{\|{\tt p}(n)\|^r} \sum_{i_1=1}^N \sum_{i_2=1}^N \cdots \sum_{i_r=1}^N \bigg\| \frac{\partial^{r}(f-g) }{\partial x_{i_1}\cdots \partial x_{i_{r}}} \bigg\|_\infty
\sum_{\mu=0}^r \left(\begin{array}{c} r\\ \mu\end{array} \right) M_\mu(\varphi) \widetilde{M}_{r-\mu}(\psi).
\end{eqnarray*}
For $J_2^2,$ since 
\begin{eqnarray*}
&& |u_{i_1}- x_{i_1}|^{r+1} \leq \|{\tt u} - {\tt x}\|^{r+1} \leq \frac{1}{\|{\tt p}(n)\|^{r+1}} (\|{\tt p}(n){\tt u}- {\tt k}\| +  \|{\tt k}- {\tt p}(n) {\tt x}\|)^{r+1}\\ &=& \frac{1}{\|{\tt p}(n)\|^{r+1}} \sum_{\mu=0}^{r+1} \left(\begin{array}{c} r+1\\ \mu\end{array} \right)\|{\tt p}(n){\tt u}- {\tt k}\|^{r+1-\mu}\|{\tt k}- {\tt p}(n) {\tt x}\|^\mu 
 \end{eqnarray*}
 again we have
\begin{eqnarray*}
&& J_2^2 \leq \frac{N}{(r+1)! \|{\tt p}(n)\|^{r+1}} \sum_{i_1=1}^N \sum_{i_2=1}^N \cdots \sum_{i_r=1}^N \max_{i_{r+1}= 1, \cdots, N} \bigg\| \frac{\partial^{r+1} g}{\partial x_{i_1}\cdots \partial x_{i_{r+1}}} \bigg\|_\infty\\ &\times& \sum_{\mu=0}^{r+1} \left(\begin{array}{c} r+1\\ \mu\end{array} \right) \sum_{{\tt k} \in \mathbb{Z}^N} |\varphi ({\tt p}(n){\tt x} - {\tt k})| \|{\tt k}- {\tt p}(n) {\tt x}\|^\mu 
\int_{\mathbb{R}^N} |\psi( {\tt t})| \| {\tt t}\|^{r+1-\mu} d{\tt t} \bigg]\\ &\leq&
\frac{N^{r+2}}{(r+1)! \|{\tt p}(n)\|^{r+1}} \sum_{i_1=1}^N \sum_{i_2=1}^N \cdots \sum_{i_r=1}^N \max_{i_{r+1}= 1, \cdots, N} \bigg\| \frac{\partial^{r+1} g}{\partial x_{i_1}\cdots \partial x_{i_{r+1}}} \bigg\|_\infty  \\&\times& \sum_{\mu=0}^{r+1} \left(\begin{array}{c} r+1\\ \mu\end{array} \right) M_\mu(\varphi)\widetilde{M}_{r+1-\mu}(\psi).
 \end{eqnarray*}
 Then as before, taking the infimum over $g \in C^{(r+1)}$ we obtain
$$J_2 \leq \frac{2}{r!}\frac{N^rA_r}{\|{\tt p}(n)\|^r} \sum_{i_1=1}^N \sum_{i_2=1}^N \cdots \sum_{i_r=1}^N 
 K \bigg( \frac{N^{2}}{2(r+1) \|{\tt p}(n)\|} \frac{B_r}{A_r}, \frac{\partial^{r}f }{\partial x_{i_1}\cdots \partial x_{i_{r}}} \bigg),$$
where $A_r= \sum_{\mu=0}^r \left(\begin{array}{c} r\\ \mu\end{array} \right) M_\mu(\varphi) \widetilde{M}_{r-\mu}(\psi)$ and $B_r=\sum_{\mu=0}^{r+1} \left(\begin{array}{c} r+1\\ \mu\end{array} \right) M_\mu(\varphi)\widetilde{M}_{r+1-\mu}(\psi). \Box$ 
 \vskip0,4cm
As a consequence of Theorem \ref{kfunctional} for functions $f\in C^2$ or $f\in C^{(r+1)}$ respectively, we have the following direct estimate 
\begin{Corollary}\label{voronov2}
We have
\begin{description}
\item[i)]
Let $f \in C^2$~ and let ${\tt x}\in \mathbb{R}^N.$~ Under the assumptions of Theorem \ref{kfunctional} i), there holds 
\begin{eqnarray*}
&&\bigg|n((S_{{\tt p}(n)}^{\varphi, \psi}f)({\tt x}) - f({\tt x})) - \sum_{|{\tt h}| = 1} {\tt a}^{\tt h} D^{{\tt h}}f({\tt x}) 
\sum_{{\tt s}= {\tt 0}}^{\tt h} \left(\begin{array}{c} {\tt h}\\ {\tt s}\end{array} \right) m^{|{\tt h-s}|}_{{\tt h-s}}(\varphi)  \widetilde{m}^{|{\tt s}|}_{\tt s}(\psi)\bigg|\\ &\leq& \sum_{|{\tt h}| = 1}|D^{{\tt h}}f({\tt x})|\left| \bigg(\frac{n}{({\tt p}(n))^{{\tt h}}} -{\tt a}^{\tt h}\bigg)\right| \bigg|\sum_{{\tt s}= {\tt 0}}^{\tt h} \left(\begin{array}{c} {\tt h}\\ {\tt s}\end{array} \right) m^{|{\tt h-s}|}_{{\tt h-s}}(\varphi)  \widetilde{m}^{|{\tt s}|}_{\tt s}(\psi)\bigg|\\ &+&
\frac{N^3 n}{2\|{\tt p}(n)\|^2} B_1 \sum_{i =1}^N \max_{j=1,\cdots,N}\bigg \|\frac{\partial^2 f}{\partial x_i\partial x_j} \bigg\|_\infty,
\end{eqnarray*}

\item[ii)]
Let $f \in C^{(r+1)}$~ and ${\tt x}\in \mathbb{R}^N.$ Under the assumptions of Theorem 2 ii) there holds 
\begin{eqnarray*}
&&\bigg|n^r((S_{{\tt p}(n)}^{\varphi, \psi}f)({\tt x}) - f({\tt x})) - \sum_{|{\tt h}| = r} {\tt a}^{{\tt h}}\frac{D^{{\tt h}}f({\tt x})}{{\tt h}!} \sum_{{\tt s}= {\tt 0}}^{\tt h} \left(\begin{array}{c} {\tt h}\\ {\tt s}\end{array} \right) m^{|{\tt h-s}|}_{{\tt h-s}}(\varphi)  \widetilde{m}^{|{\tt s}|}_{\tt s}(\psi)\bigg|\\ &\leq& 
\sum_{|{\tt h}| = r} \bigg |\frac{D^{{\tt h}}f({\tt x})}{\tt h!}\bigg| \left| \bigg(\frac{n^r}{({\tt p}(n))^{{\tt h}}} -{\tt a}^{\tt h}\bigg)\right| \bigg|\sum_{{\tt s}= {\tt 0}}^{\tt h} \left(\begin{array}{c} {\tt h}\\ {\tt s}\end{array} \right) m^{|{\tt h-s}|}_{{\tt h-s}}(\varphi)  \widetilde{m}^{|{\tt s}|}_{\tt s}(\psi)\bigg|\\ &+&
\frac{N^{r+2} n^r B_r}{(r+1)! \|{\tt p}(n)\|^{r+1}} \sum_{i_1=1}^N \sum_{i_2=1}^N \cdots \sum_{i_r=1}^N 
 \max_{i_{r+1}= 1, \cdots, N} \bigg\| \frac{\partial^{r+1} f}{\partial x_{i_1}\cdots \partial x_{i_{r+1}}} \bigg\|_\infty.
\end{eqnarray*}

\end{description}
\end{Corollary}
{\bf Proof.} 
For part i), if $f\in C^2$ we can write 
\begin{eqnarray*}
&&K\bigg(\frac{B_1}{A_1}\frac{N^2}{4\|{\tt p}(n)\|}, \frac{\partial f}{\partial x_i}\bigg) =
\inf_{g\in {C}^2}\bigg( \bigg\|\frac{\partial (f-g)}{\partial x_i}\bigg\|_\infty + \frac{B_1}{A_1}\frac{N^2}{4\|{\tt p}(n)\|}  \max_{j=1,\cdots,N}\bigg \|\frac{\partial^2 g}{\partial x_i\partial x_j} \bigg\|_\infty\bigg)\\ &\leq&
\frac{B_1}{A_1}\frac{N^2}{4\|{\tt p}(n)\|}  \max_{j=1,\cdots,N}\bigg \|\frac{\partial^2 f}{\partial x_i\partial x_j} \bigg\|_\infty
\end{eqnarray*}
by choosing $g= f.$ Therefore the assertion easily follows. In an analogous way we obtain part ii). $\Box$

\section{Examples}

In this section we apply the previous theory to some particular examples of kernels $\varphi$ and $\psi.$
We limit the study to the Voronovskaja formula, their quantitative versions are obtained in a similar way.

Here for a function $g\in L^1(\mathbb{R}^N)$ the Fourier transform of $g$ is defined as
$$\widehat{g}({\tt v}) = \int_{\mathbb{R}^N} g({\tt x}) e^{-i{\tt x}\cdot {\tt v}} d{\tt x}, \quad {\tt v}\in \mathbb{R}^N.$$
The following result will be useful (see Lemma 3.2 in \cite{BFS}, see also \cite{BM2})
\newtheorem{Proposition}{Proposition}
\begin{Proposition}\label{moment}
Let $\varphi\in C^{0}$ be such that $M_r(\varphi)<+\infty,$ for some $r\in \mathbb{N}_0$ and let ${\tt h}\in \mathbb{N}_0^N$ be fixed with $|{\tt h}|\leq r.$
The following two assertions are equivalent, for $c\in \mathbb{R}:$
\begin{description}
\item[(i)] 
$$ \sum_{{\tt k} \in \mathbb{Z}^N}\varphi ({\tt u} - {\tt k}) ({\tt k} - {\tt u})^{\tt h}= c \quad \mbox{a.e. in}\quad \mathbb{R}^N,$$

\item[(ii)] 
$$ D^{\tt h} \widehat{\varphi}(2 {\tt k}\pi)= 
 \left\{ \begin{array}{ll} (-i)^{|\tt h|} c,  & {\tt k}= {\tt 0} \\
 0,  & {\tt k} \in \mathbb{Z}^N\setminus \{\tt 0\}. 
 \end{array} \right.
 $$
\end{description}
\end{Proposition}
As a consequence, we have $m_{\tt h}^{|\tt h|}(\varphi) = (-1)^{|\tt h|}c.$
 
\begin{description}
\item[I] Let us consider as function $\varphi$ the Bochner-Riesz kernel defined by (see e.g. \cite{SW})
$$\varphi({\tt x}) \equiv b^\gamma({\tt x}) = \frac{2^\gamma}{(\sqrt{2\pi})^N}\Gamma(\gamma+1)(|{\tt x}|)^{-N/2 -\gamma} J_{(N/2)+ \gamma}(|{\tt x}|)$$
for $\gamma>0,$ where $J_\lambda$ is the Bessel function of order $\lambda.$
It is well known that
$$ \widehat{b^\gamma}({\tt v}) = \left\{ \begin{array}{ll} (1- |{\tt v}|^2)^\gamma,  & |{\tt v}|\leq 1 \\
 0,  & |{\tt v}|>1 . 
 \end{array} \right.$$
Using Proposition \ref{moment} we obtain, for every ${\tt u}\in \mathbb{R}^N,$ 
$$m^0_{{\tt 0}}(\varphi)= \sum_{{\tt k} \in \mathbb{Z}^N}  b^\gamma({\tt u} - {\tt k})= \widehat{b^\gamma}({\tt 0})=1$$
and for $i, j=1, \cdots, N$
$$m^1_{{\tt e}_i}(b^\gamma)=0, \quad m^2_{2{\tt e}_i}(b^\gamma)=2\gamma, \quad m^2_{{\tt e}_i+ {\tt e}_j}(b^\gamma)=0,\quad i\neq j.$$
Moreover for $\gamma >5/2$ we have that (see \cite{BM2}) $M_2(b^\gamma)<+\infty$ and
$$\lim_{w\rightarrow +\infty} \sum_{\|{\tt u} - {\tt k} \|>w} |b^\gamma ({\tt u} - {\tt k})\|{\tt u} - {\tt k}\|^2=0,$$
uniformly with respect to ${\tt u}\in \mathbb{R}^N.$

For what concerns the kernel $\psi,$ let $B_2$ the one-dimensional central B-splines of order 2 defined by
$$B_2(x) = (1- |x|)\chi_{]-1,1[}(x)\quad x\in \mathbb{R}$$
and let us consider the kernel 
$$\psi({\tt x})= B_2(x_1)B_2(x_2)\cdots B_2(x_N).$$
Since $\int_{\mathbb{R}} B_2(x)dx=1$, $\int_{\mathbb{R}} x^2B_2(x)dx=\frac{1}{6}$
and $\int_{\mathbb{R}} x^k B_2(x)dx=0$ for every $k$ odd,
we have 
$$\widetilde{m}^0_{{\tt 0}}(\psi) = \int_{\mathbb{R}^N}\psi({\tt x}) d{\tt x} =1$$
and for $i, j=1, \cdots, N$
$$\widetilde{m}^1_{{\tt e}_i}(\psi)=0, \quad \widetilde{m}^2_{2{\tt e}_i}(\psi)=\frac{1}{6}, \quad \widetilde{m}^2_{{\tt e}_i+ {\tt e}_j}(\psi)=0,\quad i\neq j.$$
Since $\psi$ has compact support all the assumptions i), ii) and iii) for the kernels are satisfied with $r=2.$
Thus we can apply Corollary \ref{voronov} ii) with $r=2,$ (see Remark 2), obtaining
\begin{eqnarray*}
&&\lim_{n \rightarrow + \infty}n^2 [(S_{{\tt p}(n)}^{\varphi, \psi}f)({\tt x}) - f({\tt x})] = (\gamma + \frac{1}{12}) 
\sum_{k=1}^N a_k^2\frac{\partial^2 f({\tt x})}{\partial x_k^2}.
\end{eqnarray*}

\item[II] Let us consider the kernel $\varphi({\tt x})= \psi({\tt x})= B_2(x_1)B_2(x_2)\cdots B_2(x_N).$
In this case the moment of the kernel $\varphi$ are (see \cite{BFM})
$$ m^0_{{\tt 0}}(\varphi)=1, \quad  m^1_{{\tt e}_i}(\varphi)=0, \quad  m^2_{2{\tt e}_i}(\varphi)=\frac{1}{6}, \quad m^2_{{\tt e}_i+ {\tt e}}(\varphi)= 0$$
and also in this case we can apply Corollary \ref{voronov} ii) with $r=2$ (see Remark 2) obtaining
\begin{eqnarray*}
&&\lim_{n \rightarrow + \infty}n^2 [(S_{{\tt p}(n)}^{\varphi, \psi}f)({\tt x}) - f({\tt x})] = \frac{1}{6}
\sum_{k=1}^N a_k^2 \frac{\partial^2 f({\tt x})}{\partial x_k^2}.
\end{eqnarray*}

\item[III] Let us put $\varphi({\tt x})=B_2(x_1)B_2(x_2)\cdots B_2(x_N)$ and $\psi({\tt x}) = F(x_1)\cdots F(x_N)$ where 
the function $F: \mathbb{R}\rightarrow \mathbb{R}$  is defined by
$$ F(x)= \left\{ \begin{array}{ll} x+1,  & -1\leq x \leq 0 \\
 \frac{1}{2} e^{-x},  & x>0\\
 0, & \mbox{otherwise} . 
 \end{array} \right.$$
 
Since $\int_{\mathbb{R}} F(x)dx=1$, $\int_{\mathbb{R}} x F(x)dx=\frac{1}{3}$, we have for $i=1, \cdots, N$
$$\widetilde{m}^0_{{\tt 0}}(\psi) = \int_{\mathbb{R}^N}\psi({\tt x}) d{\tt x} =1, \quad \widetilde{m}^1_{{\tt e}_i}(\psi)=\frac{1}{3},$$ 
so we can apply Corollary \ref{voronov} i) with $r=1$ obtaining
\begin{eqnarray*}
&&\lim_{n \rightarrow + \infty}n [(S_{{\tt p}(n)}^{\varphi, \psi}f)({\tt x}) - f({\tt x})] = \frac{1}{3}\sum_{k=1}^N a_k\frac{\partial f}{\partial x_k}({\tt x}).
\end{eqnarray*}

\end{description}
\vskip0,3cm
\noindent
{\bf Aknowledgments}. The authors have been partially supported by the Gruppo Nazionale Analisi Matematica, Probabilit\'a e Applicazioni (GNAMPA) of the Istituto Nazionale di Alta Matematica (INdAM) and by the Department of Mathematics and Computer Science of University of Perugia.


\end{document}